\documentclass[12pt]{article}
\usepackage[english]{babel}
\usepackage{amssymb,amsmath,graphics}
\parskip\medskipamount
\newtheorem{theorem}{Theorem}[section]
\newtheorem{remark}[theorem]{Remark}
\newtheorem{prop}[theorem]{Proposition}
\newtheorem{lemma}[theorem]{Lemma}
\newtheorem{mlemma}[theorem]{Main Lemma}
\newtheorem{cor}[theorem]{Corollary}
\def\<{\langle}
\def\>{\rangle}
\newcommand{\proof}{\emph{Proof.~}}
\newcommand{\cC}{\mathcal{C}}

\newcommand{\cA}{\mathcal{A}}

\newcommand{\cM}{\mathcal{M}}

\newcommand{\st}{~|~}

\def\qed{{\hfill\hphantom{.}\nobreak\hfill$\Box$}}
\newcommand{\proj}{\mathrm{proj}}
\newcommand{\supp}{\mathrm{supp}}
\newcommand{\id}{\mathrm{id}}
\renewcommand{\d}{\mathrm{d}}

\begin{document}

\author{Peter Abramenko \and Hendrik Van Maldeghem\thanks{The second author is partly supported by a
Research Grant of the Fund for Scientific Research - Flanders (FWO - Vlaanderen)}}
\title{\bf Intersections of Apartments}
\date{}
\maketitle

\begin{abstract}We show that, if a building is endowed with its complete system of apartments, and if each
panel is contained in at least four chambers, then the intersection of two apartments can be any convex
subcomplex contained in an apartment. This combinatorial result is particularly interesting for lower dimensional convex
subcomplexes of apartments, where we definitely need the assumption on the four chambers per panel in the building.
The corresponding statement is not true anymore for arbitrary systems of apartments, and
counter-examples for infinite convex subcomplexes exist for any type of buildings. However, when we restrict to
finite convex subcomplexes, the above remains true for arbitrary systems of apartments if and only if every
finite subset of chambers of the standard Coxeter complex is contained in the convex hull of two chambers.
\end{abstract}

\section{Introduction}
\emph{Buildings} were defined by Jacques Tits in the sixties as the natural geometric structures related to
groups of Lie type -- providing natural permutation modules for these groups. In this setting, only the
spherical buildings are important, and a full classification for irreducible rank at least three exists, see the seminal
monograph \cite{Tit:74}. However, it became clear later that also non-spherical buildings (e.g.~affine buildings) play a
very important role in group theory, and we will not make any restrictions concerning the type of the building in this paper.
In the original definition of buildings apartments play a crucial role. In fact, the two main
axioms of buildings in \cite{Tit:74} are about apartments. One axiom states that every pair of chambers is contained
in an apartment, and the other that two apartments can be isomorphically mapped onto each other by an
isomorphism fixing two chosen simplices in the intersection and all their faces. Naturally, one can ask which convex
subcomplexes of apartments can be realized as intersections of two apartments. We provide a complete answer to this
question for thick buildings in the present paper.  Note that this result goes a
little bit in the opposite direction of the main use of apartments. Indeed, usually one reduces a problem concerning
buildings to a problem in an apartment | a Coxeter complex. Here, we reduce potential problems involving the
intersection of roots in an apartment (which is by definition a convex subcomplex of that
apartment) to the intersection of just two apartments, but in the whole building.

In the eighties, a local approach to buildings \cite{Tit:81}, again initiated by Jacques Tits, led the
latter to introduce buildings without mentioning apartments, but using chamber systems. But
around the same time, or even a little bit later, Jacques Tits completed the classification of affine buildings
of irreducible rank at least four, and there the apartments play a prominent role \cite{Tit:86}. In fact, they were so
important that the class of objects to which the classification could be applied (this class is slightly larger
than that of the affine buildings) was called \emph{systems of apartments} by Tits. Later, in the nineties, yet another
definition of buildings was introduced, again by Jacques Tits, this time for the benefit of
\emph{twin buildings}, which are natural generalizations of spherical buildings, and which are the natural
geometric structures for Kac-Moody groups \cite{Tit:92} over fields. That new definition interprets buildings as $W$-metric spaces,
a point of view which will also be important in parts of the present paper (see e.g.~Lemma~\ref{lemma1}). It again does not explicitly mention
apartments which, however, remain to be an important notion also in this approach . For instance, together with twin
buildings, the notion of a  \emph{generalized Moufang building} became important.
This notion, which was first introduced for spherical buildings in the appendix of \cite{Tit:74}, uses in a
crucial way so-called \emph{roots} in a building, which are nothing other than \ldots half-apartments.

So it is clear that the apartments of a building are of vital importance in the whole theory. Although one
sometimes removes them from the definition, in order to make things simpler, they have remained crucial in
the results and the theory. Therefore, the question of what the intersection of two apartments can be is a
fundamental one for all types of buildings. We answer this question in the present paper. The result is
neither surprising nor difficult to prove if this intersection contains chambers. However, the problem becomes
much more difficult in the case of lower dimensional convex subcomplexes of apartments. First of all, it is clear
that then not every convex subcomplex is an intersection of two apartments if the building has only three chambers
per panel. It is quite surprising that this is in fact the only obstacle, i.e., in a building with at least four
chambers per panel, \emph{any} convex subcomplex of an apartment is an intersection of two apartments.  In this
context `apartment' means `member of the complete system of apartments'. If one also considers arbitrary apartment
systems, then the answer becomes more involved, and we have also dealt with that situation. The simple case
of two adjacent chambers | which is for every building and every apartment system a convex set contained in an
apartment | already has a
nice application, as we will show.

The paper is organized as follows. In Section~\ref{MR} we introduce some basic notions and fix some notation. We
also state our main results and an application. In the rest of the paper, we then proceed to prove these
results. Along the way, we have phrased some lemmas and propositions slightly more generally then necessary for
our purposes, because we think that these results can be of independent use. For instance,
Section~\ref{constructions} deals with constructing apartments under various conditions | obviously important
in this paper, but potentially also useful in many other situations. In Section~\ref{complete}, we consider
the case of a complete system of apartments, while Section~\ref{incomplete} looks at the other situations.

\section{Preliminaries and Main Results}\label{MR}

We assume that the reader is familiar with basic concepts from
building theory, as presented in \cite{Abr-Bro:08}, \cite{Ron:89}
and \cite{Tit:74}. In order to fix some notation and avoid
confusion, we will repeat some of these concepts here. We mainly
take the simplicial point of view. So a building $\Delta$ will be a
chamber complex which admits a system of apartments satisfying the
usual axioms (see e.g.~\cite{Tit:74}, Chapter 3, or
\cite{Abr-Bro:08}, Section~4.1). However, on the one hand, we shall
\emph{not} always assume that $\Delta$ is thick, but state this
condition explicitly whenever needed. On the other hand, we
\emph{do} assume that $\Delta$ is finite-dimensional. As usual, the
top-dimensional simplices of $\Delta$ are called \emph{chambers},
the simplices of codimension 1 are called \emph{panels}, and
$\Delta$ is called thick if each panel is contained in at least
three chambers. We denote by $\cC=\mathrm{Ch}(\Delta)$ the set of
chambers of $\Delta$, and, for any chamber subcomplex $\kappa$ of $\Delta$,
by $\mathrm{Ch}(\kappa)$ the set of chambers of $\kappa$. Two
chambers are called \emph{adjacent} if they are different (note that
some authors do not require this) and share a common panel. A
\emph{gallery (of length $n$)} is a sequence $(C_0,C_1,\ldots,C_n)$
of $n+1$ chambers in which $C_{i-1}$ and $C_i$ are adjacent, for
all $1\leq i\leq n$. Such a gallery is called \emph{minimal} if
there is no gallery of length $<n$ starting in $C_0$ and ending in
$C_n$. In this case we also say that $n$ is the \emph{(gallery)
distance} between $C_0$ and $C_n$. If two chambers $C,D\in\cC$ are at distance $n$, we write $\d(C,D)=n$.

A chamber subcomplex $\kappa$ of $\Delta$ is called \emph{convex} if
for any two chambers $C,D$ of $\kappa$, every gallery
$(C_0,C_1,\ldots,C_n)$ in $\Delta$ with $C_0=
C$, $C_n=D$ and $n=\d(C,D)$ is contained in $\kappa$. 
An arbitrary subcomplex $\kappa$ of $\Delta$ is called
\emph{convex} if it is the intersection of convex chamber
subcomplexes of $\Delta$. This is the definition of convexity given
by Tits in \cite{Tit:74}, Section~1.5, but differs from the one given in
\cite{Abr-Bro:08}, Section 4.11. However, if $\kappa$ contains a
chamber or is contained in an apartment of $\Delta$, these two
definitions agree, and in this paper we will only consider
subcomplexes of apartments. For any set $\mathcal{S}$ of simplices
of $\Delta$, the \emph{convex closure} or \emph{convex hull} of
$\mathcal{S}$ is defined as the intersection of all convex subcomplexes
containing $\mathcal{S}$.

Recall that all apartments of $\Delta$ are Coxeter complexes which
are isomorphic to each other. If they are isomorphic to the standard
Coxeter complex $\Sigma(W,S)$ (as defined in \cite{Abr-Bro:08},
Section 3.1), where $(W,S)$ is a Coxeter system, we say that
$\Delta$ is of type $(W,S)$. Each building admits a maximal
apartment system $\widetilde{\cA}$ (see for instance \cite{Abr-Bro:08}, Section~4.5), 
which is called the
\emph{complete} system of apartments of $\Delta$. A subset $\cA$ of
$\widetilde{\cA}$ is a system of apartments of $\Delta$ if and only
if any two chambers of $\Delta$ are contained in some element of
$\cA$. If $\Delta$ is of type $(W,S)$, the elements of
$\widetilde{\cA}$ are precisely the simplicial subcomplexes of
$\Delta$ which are isomorphic to $\Sigma(W,S)$ (see
\cite{Abr-Bro:08}, Proposition~4.59). Whenever we talk about
apartments of $\Delta$ without further specification, we mean
elements of the complete system $\widetilde{\cA}$ of apartments of
$\Delta$.

We fix an apartment $\Sigma$ of a building $\Delta$. As a preparation towards our First Main Result, we shall prove in Section~\ref{complete}:

\begin{prop}\label{abr}
If $\Delta$ is thick and $\kappa$ is a convex chamber subcomplex of $\Sigma$, then there exists an apartment $\Sigma'$ of
$\Delta$ such that $\kappa=\Sigma\cap\Sigma'$.
\end{prop}

\begin{remark}\rm
Assume there exists a panel $P$ in $\Sigma$ which is contained in exactly $3$ chambers of $\Delta$. Then it is
obvious that the wall $M$ of $\Sigma$ containing $P$ cannot be written as $M=\Sigma\cap\Sigma'$ for some
apartment $\Sigma'$ since $\Sigma\cap\Sigma'$ has to contain at least one of the three chambers
containing $P$.\end{remark}

So it is clear that Proposition~\ref{abr} cannot be true in general for (lower dimensional) convex subcomplexes
of $\Sigma$. Amazingly enough, it turns out to be true if we only add the (obviously necessary) condition that
each panel is contained in at least $4$ chambers of $\Delta$.

\begin{theorem}[First Main Result]\label{mr1}
Suppose every panel of the building $\Delta$ is contained in at least $4$ chambers. Then for every convex
subcomplex $\kappa$ of $\Sigma$, there exists an apartment $\Sigma'$ of $\Delta$ such that
$\kappa=\Sigma\cap\Sigma'$.\end{theorem}

Now a natural question arises: What can we say if the building $\Delta$ is \emph{not} endowed with its
complete set of apartments? We will prove that, whenever a convex subcomplex of an apartment has
infinitely many vertices, then Theorem~\ref{mr1} does not hold in general (when the system of apartments is
not complete). In fact, in this case, such a convex subcomplex might not even be the intersection of
\emph{all} apartments in which it is contained. Hence we are lead to consider only the case where our
convex subcomplex has a finite number of simplices. Then we will show that every such convex subcomplex
(contained in an apartment of an arbitrarily given apartment system $\cA$ of $\Delta$) is the intersection of two 
apartments of $\cA$ if and only if the building is a direct product
of spherical and affine buildings. Here is the precise result:
\begin{theorem}[Second Main Result]\label{mr2}
Let $\Delta$ be a thick building of type $(W,S)$. Then the following are equivalent.
\begin{itemize}
\item[\rm(I)] For any apartment system $\cA$ of $\Delta$, the convex closure of two given chambers $C$ and $D$
is the intersection of two apartments in $\cA$.

\item[\rm(II)] For any apartment system $\cA$ of $\Delta$, every convex finite set $\kappa$ of chambers
contained in some member of $\cA$ is the intersection of two apartments in $\cA$.

\item[\rm(III)] For every finite subset $F$ of chambers of the Coxeter complex $\Sigma(W,S)$, there exist
two chambers $C,D\in\Sigma(W,S)$ such that $F$ is contained in the convex closure of $C$ and $D$.
\end{itemize}

If every panel of $\Delta$ is contained in at least four chambers, and Condition $(III)$ is satisfied, then
every finite convex subcomplex contained in an apartment is the intersection of two apartments (with respect to any
apartment system $\cA$).
\end{theorem}

Now it is proved by Caprace in \cite{Cap:05} that Assertion~(III) is equivalent to $\Sigma(W,S)$ having
only spherical and affine components.

We shall also prove the following related result for pairs of adjacent chambers.

\begin{prop}\label{propmain}Let $\Delta$ be a thick building endowed with an arbitrary system $\cA$ of apartments.
Then the set of two adjacent chambers and all their faces is always the intersection of all apartments in $\cA$ in which they are both contained.
\end{prop}

As an application of this proposition we shall prove in Section~\ref{incomplete}:

\begin{prop}
Every map between the chamber sets of a thick building $\Delta$ and an arbitrary building $\Delta'$ (endowed
with  arbitrary systems of apartments $\cA$ and $\cA'$, respectively) which bijectively maps the set of
chambers of any apartment $\Sigma\in\cA$ to the set of chambers of some apartment $\Sigma'\in\cA'$ is injective and
preserves adjacency of chambers. Hence it induces a simplicial isomorphism of $\Delta$ onto a thick subbuilding
of $\Delta'$ in case $\Delta$ is $2$-spherical.
\end{prop}

\section{Constructing apartments}\label{constructions}

In the course of the proof of our main results, we shall make use of some well-known facts, and some lemmas
that might be of independent interest. We collect these assertions in this section. Most of them are concerned
with the construction of apartments containing certain sets of simplices and satisfying certain conditions.

We start by repeating some more notions from the simplicial theory of buildings. If $\Sigma$ is a Coxeter complex, a \emph{root} of $\Sigma$ is the image of a folding of $\Sigma$ (see \cite{Tit:74}, Chapter 2, or \cite{Abr-Bro:08}, Section 3.4). It is well known that a subcomplex $\kappa$ of $\Sigma$ is convex if and only if it is an intersection of roots of $\Sigma$. If $\alpha$ is a root of $\Sigma$, the subcomplex of $\Sigma$ generated by all chambers of $\Sigma$ which are not in $\alpha$ is again a root of $\Sigma$, denoted by $-\alpha$ (if $\Sigma$ is understood) and called the root \emph{opposite} $\alpha$ in $\Sigma$. The subcomplex $\alpha\cap(-\alpha)$ is called the \emph{boundary} or \emph{wall} of $\alpha$ (or of $-\alpha$); it is the subcomplex of $\Sigma$ generated by all panels in $\Sigma$ which are contained in precisely one chamber of $\alpha$. If $\Delta$ is a building, a \emph{root} or \emph{wall} in $\Delta$ is a root or wall in one of its apartments.

In a Coxeter complex $\Sigma$, the \emph{support} $\supp A$ of a simplex $A$ is the intersection of all walls of $\Sigma$ containing $A$. We repeat the following result from \cite{Abr:94}:

\begin{prop}\label{fact1}
If $\kappa$ is a convex subcomplex of a Coxeter complex $\Sigma$, $A$ is a maximal simplex of $\kappa$, and $M$ is a wall in
$\Sigma$ containing $A$, then $\kappa\subseteq M$. Hence, $\kappa\subseteq \supp A$.
\end{prop}

\proof See Proposition 1 $(ii)$ of \cite{Abr:94}.

In the following, $\Delta$ denotes a building of type $(W,S)$, and $\cC$ is the set of chambers of $\Delta$.
The following property of apartments is well known if $\Sigma\cap\Sigma'$ contains chambers, and is proved in general in \cite{Abr-Bro:08}:

\begin{prop}\label{newprop}
Let $\Sigma,\Sigma'$ be two apartments of $\Delta$. Then there exists an isomorphism
$\varphi:\Sigma\rightarrow\Sigma'$ such that $\varphi$ fixes each simplex in $\Sigma\cap\Sigma'$.
\end{prop}

\proof See Proposition 4.101 of \cite{Abr-Bro:08}

Before we continue, we remind the reader that $\Delta$
comes equipped with a
\emph{Weyl distance function} $\delta:\cC\times\cC\rightarrow W$
(see \cite{Abr-Bro:08}, Section 4.8). This function $\delta$ is an
important tool in the theory of buildings. The latter can, in fact,
be based completely on the properties of the Weyl distance function
(see \cite{Abr-Bro:08}, Chapter 5). If $\Delta$ and $\Delta'$ are
two buildings of type $(W,S)$ with chamber sets $\cC$ and $\cC'$, and with Weyl distance functions $\delta$ and $\delta'$,
respectively, then two subsets $\mathcal{M}\subseteq\cC$ and
$\mathcal{M}'\subseteq \cC'$ are called \emph{$W$-isometric} (or just
\emph{isometric}) if there exists a bijective map
$f:\mathcal{M}\rightarrow\mathcal{M}'$ such that
$\delta'(f(C),f(D))=\delta(C,D)$, for all $C,D\in\mathcal{M}$. One of
the fundamental facts in this context is that a subset
$\mathcal{M}\subseteq\cC$ is contained in $\mathrm{Ch}(\Sigma)$, for
some apartment $\Sigma$ of $\Delta$, if and only if $\mathcal{M}$
is isometric to a subset of $W=\mathrm{Ch}(\Sigma(W,S))$, see
\cite{Ron:89}, Theorem 3.6, or \cite{Abr-Bro:08}, Theorem 5.73.

This will be used in the next two propositions.

\begin{prop}\label{extension}
Given a root $\alpha$ of $\Delta$ with boundary wall $M$, a panel $P\in M$ and a chamber $C$ containing $P$ but not
contained in $\alpha$, there exists an apartment $\Sigma$ containing $\alpha\cup\{C\}$.
\end{prop}

\proof Consider any apartment $\Sigma'$ containing $\alpha$ and let $D$ be the chamber in $\Sigma'$ containing
$P$ and not contained in $\alpha$. Then it is clear that  $\mathrm{Ch}(\alpha)\cup\{C\}$ is $W$-isometric to
$\mathrm{Ch}(\alpha)\cup\{D\}$. The assertion now follows from Theorem 3.6 in \cite{Ron:89}. \qed

\begin{prop}\label{fact2}
A subset $\cM$ of $\cC=\mathrm{Ch}(\Delta)$ is contained in the set of chambers of some apartment of $\Delta$
if and only if \begin{itemize}\item[$(\star)$] $\delta(C,E)=\delta(C,D)\delta(D,E)$, for all $C,D,E\in\cM$.
\end{itemize}
\end{prop}

\proof This is contained in \cite{Abr-Bro:08}, see Exercise 5.77. We provide a proof for completeness' sake.

The necessity of Condition~$(\star)$ is clear by standard properties of apartments. To verify sufficiency,
first fix a chamber $C_0\in\cM$ and define a map $f:\cM\rightarrow W$ by $f(C)=\delta(C_0,C)$ for all
$C\in\cM$. Due to the Condition~$(\star)$, this is an isometry. Indeed, denoting the natural Weyl distance on
$W$ by $\delta_W$, we have
\begin{eqnarray*}
\delta_W(f(C),f(D))=f(C)^{-1}f(D)&=&\delta(C_0,C)^{-1}\delta(C_0,D)\\
                                 &=&\delta(C,C_0)\delta(C_0,D)=\delta(C,D).
\end{eqnarray*}

Slightly abusing notation, we also denote the induced map $\cM\rightarrow f(\cM)$ by $f$. Then $f$ is a
bijective isometry, so that $f^{-1}:f(\cM)\rightarrow \cM$ is an isometry. Now, by Theorem~3.6 of \cite{Ron:89}, $f^{-1}$ can be
extended to an isometry $g:W\rightarrow\cC$, and hence $\cM$ is contained in $g(W)$, which is the set of
chambers of an apartment of $\Delta$. \qed

We now prove some lemmas. Their proofs use projections, a concept which we briefly recall now (for details, 
we refer to \cite{Tit:74}, Chapter 3, or \cite{Abr-Bro:08}, Sections 4.9 and~5.3).
Let $A,B$ be simplices of the building $\Delta$. If A and $B$ are chambers, we denote by $\d(A,B)$ the gallery 
distance between $A$ and $B$. In general, we define $\d(A,B)=\min\{\d(C,D)~|~C,D\in\cC,
A\subseteq C\mbox{ and }B\subseteq D\}$. Then there is a unique simplex $P$ containing $A$ with the 
property that, for any chamber $C$ containing $A$, one has $\d(C,B)=\d(A,B)$ if and only if
$P\subseteq C$. This simplex $P$ is called the \emph{projection of $B$ onto $A$} and denoted by $P=\proj_AB$. 
If $B$ is a chamber, then also $\proj_AB$ is a chamber, and it has the ``gate property''
$\d(C,B)=\d(C,\proj_AB)+\d(\proj_AB,B)$, for all chambers $C$ containing $A$. Moreover, 
$\delta(C,B)=\delta(C,\proj_AB)\delta(\proj_AB,B)$ in this case (see \cite{Abr-Bro:08}, Proposition~5.34).

\begin{lemma}\label{lemma3} Let $\alpha_1,\alpha_2$ be two roots in $\Delta$ having the same wall $M$ as their
boundary. Then $\alpha_1\cup\alpha_2$ is an apartment of $\Delta$ if and only if there exists a panel $P\in M$
such that the chamber $C_1\in\alpha_1$ with $P\subseteq C_1$ is different from the chamber $C_2\in\alpha_2$
with $P\subseteq C_2$.
\end{lemma}

\proof Obviously, if $C_1=C_2$, then $\alpha_1\cup\alpha_2$ cannot be an apartment. So let us assume now that
$C_1\neq C_2$. We first claim that this implies $\alpha_1\cap\alpha_2=M$. Indeed, let $x$ be any simplex in
$\alpha_1$ with $x\notin M$. Then $\proj_Px\neq P$, because otherwise there
would be a minimal gallery between $P$ and $x$ in some apartment $\Sigma$ containing $\alpha_1$, starting with a chamber
outside $\alpha_1$, which is impossible since $\alpha_1$ is convex and every chamber of $\Sigma$ containing $x$ has
to be in $\alpha_1$. Hence $P$ is strictly contained in $\proj_Px$. But the only simplex in
$\alpha_1$ strictly containing $P$ is $C_1$, implying $\proj_Px=C_1$. By assumption, $C_1\notin \alpha_2$.
Hence $x$ cannot belong to $\alpha_2$. Similarly, $y\notin\alpha_1$ for any $y\in\alpha_2\setminus M$.
Therefore, $\alpha_1\cap\alpha_2=M$.

We now consider apartments $\Sigma_1,\Sigma_2\in\cA$ containing $\alpha_1,\alpha_2$, respectively. Denote by
$\beta_i$ the root opposite $\alpha_i$ in $\Sigma_i$, $i=1,2$. By Proposition~\ref{newprop}, there exists
a simplicial isomorphism $\varphi:\Sigma_1\rightarrow\Sigma_2$ which fixes $\Sigma_1\cap\Sigma_2$
pointwise. So in particular $\varphi$ is the identity on $M$. This implies that $\varphi(\alpha_1)$ has again
$M$ as its bounding wall, i.e., $\varphi(\alpha_1)=\alpha_2$ or $\varphi(\alpha_1)=\beta_2$. We may assume
$\varphi(\alpha_1)=\beta_2$ as otherwise we simply compose $\varphi$ with the reflection in $\Sigma_2$ about
the wall $M$, and this reflection interchanges $\alpha_2$ with $\beta_2$ (and note that this composite map
still fixes $M$ pointwise). We now define a simplicial map $\rho:\alpha_1\cup\alpha_2\rightarrow\Sigma_2$ by
$\rho_{|\alpha_2}=\id_{\alpha_2}$ and $\rho_{|\alpha_1}=\varphi_{|\alpha_1}$.

Note that $\rho$ is well defined since $\alpha_1\cap\alpha_2=M$ and $\varphi_{|M}=\id_M$. It is also clear that
$\rho$ is a simplicial morphism. Next we define $\sigma:\Sigma_2\rightarrow\alpha_1\cup\alpha_2$ by
$\sigma_{|\alpha_2}=\id_{\alpha_2}$ and $\sigma_{|\beta_2}=\varphi^{-1}_{|\beta_2}$. Again, $\sigma$ is a well
defined simplicial morphism (here we use $\alpha_2\cap\beta_2=M$), and it is obvious that $\sigma$ and $\rho$
are inverse to each other. Hence $\sigma:\Sigma_2\rightarrow\alpha_1\cup\alpha_2$ is a simplical isomorphism. It follows that
$\alpha_1\cup\alpha_2=\sigma(\Sigma_2)$ is an apartment of $\Delta$, by Proposition 4.59 in \cite{Abr-Bro:08} (note that 
the application of the latter does not require $\sigma$ to be type-preserving, although it is easy to see that it is here).\qed

We mention a consequence of this lemma, thereby slightly improving Lemma 2.2 of \cite{Abr-Mal:04}. We denote
the set of chambers containing a given panel $P$ by $\cC_P$.

\begin{cor} \label{cor} Given a panel $P$ and a wall $M$ in $\Delta$ with $P\in M$. There
exists a family $(\alpha_C)_{C\in\cC_P}$ of roots with the following properties.
\begin{itemize}
\item[$(1)$]$C\in\alpha_C$, for all $C\in\cC_P$;

\item[$(2)$] $M$ is the boundary wall of $\alpha_C$;

\item[$(3)$] $\alpha_C\cup\alpha_D$ is an apartment of $\Delta$,
for all $C,D\in\cC_P$ with $C\neq D$.
\end{itemize}
\end{cor}

\proof Let $\Sigma_0$ be any apartment containing $M$, and let $C_0,D_0\in\cC_P$ be the two chambers in
$\Sigma_0$ containing $P$, with $C_0\neq D_0$. Let $\alpha_{C_0},\alpha_{D_0}$ be the roots in $\Sigma_0$ which
are bounded by the wall $M$ and contain the chambers $C_0,D_0$, respectively.

Now, for any $C\in\cC_P\setminus\{C_0,D_0\}$, there exists an apartment $\Sigma_C$ containing
$\alpha_{C_0}\cup\{C\}$ (direct by Proposition~\ref{extension}). Let $\alpha_C$ be the root in
$\Sigma_C$ with bounding wall $M$ (note that $M\subseteq\alpha_{C_0}\subseteq\Sigma_C$) and containing $C$.
Given two distinct elements $C,D\in\cC_P$, the two roots $\alpha_C$ and $\alpha_D$ have the wall $M$ in common,
but they do not share the same chambers $C$ and $D$, respectively, through $P\in M$. So we can apply
Lemma~\ref{lemma3} and infer that $\alpha_C\cup\alpha_D$ is an apartment of $\Delta$. \qed

The following technical lemma will be needed in the proof of the Main Lemma in Section~\ref{complete}. Before we state it, 
we recall the definition of links and stars. If $\Omega$ is an arbitrary simplicial complex and $A$ is a simplex in $\Omega$, 
we define the \emph{link of $A$ in $\Omega$} as the subcomplex $$\mathrm{Lk}_\Omega(A)=\{B\in\Omega~|~A\cup B \mbox{ is a 
simplex in }\Omega\mbox{ and }A\cap B=\emptyset\},$$ and the \emph{star of $A$ in $\Omega$} as 
$$\mathrm{St}_{\Omega}(A)=\{B\in\Omega~|~A\cup B\mbox{ is a simplex in }\Omega\}.$$ Recall that for $A\in\Delta$, $\mathrm{Lk}_\Delta(A)$ 
is again a building (with apartments $\mathrm{Lk}_\Sigma(A)$ for all apartments $\Sigma$ of $\Delta$ containing $A$).

\begin{lemma}\label{lemma1}
Let $M$ be a wall in $\Delta$, let $A$ be a simplex of $M$ and let $L_A:=\mathrm{Lk}_\Delta(A)$ be its link in
$\Delta$. Then for every apartment $\Sigma_A$ of the building $L_A$ with $M\cap L_A\subseteq \Sigma_A$, there
exists an apartment $\Sigma$ of $\Delta$ such that $\Sigma_A=\Sigma\cap L_A$ and $M\subseteq\Sigma$.
\end{lemma}

\proof In this proof, we will have to combine the simplicial view with the $W$-metric approach to buildings.
This requires some additional notation.

We set $\cC_A=\{C\in\cC\st A\subseteq C\}$. If $J\subseteq S$ denotes the cotype of $A$, then $\cC_A$ is a
\emph{$J$-residue} in $\cC$, and $(\cC_A,\delta_{|\cC_A\times\cC_A})$ is a $W$-metric building of type
$(W_J,J)$ which corresponds to the simplicial building $L_A=\mathrm{Lk}_\Delta(A)$. To the apartment $\Sigma_A$
of $L_A$ corresponds the apartment (set of chambers) $\cC'_A:=\{A\cup B\st B\mbox{ is a chamber of }\Sigma_A\}$
of $\cC_A$. Note that, therefore, $\cC'_A$ satisfies Condition~$(\star)$ of Proposition~\ref{fact2}.

We now embark on the proof of the lemma. We first choose a (simplicial) root $\widehat{\alpha}$ in $\Delta$
which has $M$ as boundary and set $\alpha:=\mathrm{Ch}(\widehat{\alpha})$. We then choose a panel $P\in M$ with
$A\subseteq P$ and $P\setminus A\in \Sigma_A$ ( recall that $M\cap L_A\subseteq \Sigma_A$), and denote by
$C_1,C_2$ the two chambers in $\cC'_A$ which contain $P$. Set $$\alpha_i:=\{C\in\cC'_A\st
\mathrm{d}(C,C_i)<\mathrm{d}(C,C_{3-i})\},$$ for $i=1,2$. These sets $\alpha_1$ and $\alpha_2$ are the sets of
chambers corresponding to the two roots in $\Sigma_A$ with boundary wall $M\cap L_A$ (which is contained in
$\Sigma_A$ by assumption).

STEP I. We first prove the assertion in the special case that $\alpha_1\subseteq\alpha$. In this case we will show
that there exists an apartment $\Sigma$ of $\Delta$ with $\alpha\cup\alpha_2\subseteq\mathrm{Ch}(\Sigma)$.

We want to verify Condition~$(\star)$ for $\alpha\cup\alpha_2$. To that aim, we choose three chambers
$C,D,E\in\alpha\cup\alpha_2$ and distinguish some (non-disjoint) cases.

\begin{itemize}
\item[(1)] $C,D,E\in\alpha$. Condition~$(\star)$ is satisfied since $\alpha$ is part of an apartment.

\item[(2)] $C,D,E\in\alpha_1\cup\alpha_2$. As remarked above, $\alpha_1\cup\alpha_2=\cC_A'$ satisfies
Condition~$(\star)$.

\item[(3)] $C\in\alpha$; \ $D,E\in\alpha_1\cup\alpha_2$. First note that the simplicial roots in $L_A$
corresponding to $\cC_A\cap\alpha$ and $\alpha_1$ have the same wall $M\cap L_A$. Since $\alpha_1\subseteq
\cC_A\cap\alpha$, they must be equal. Now set $C':=\proj_AC\in\cC_A\cap\alpha=\alpha_1$. Then, using the
standard property of the projection mapping that $\delta(C,C')\delta(C',X)=\delta(C,X)$, for all chambers
$X\in\cC_A$ (see Proposition 5.34(2) in \cite{Abr-Bro:08}), we obtain

$$\delta(C,D)\delta(D,E)=\delta(C,C')\delta(C',D)\delta(D,E)=\delta(C,C')\delta(C',E)=\delta(C,E),$$

where we also used the fact that $\cC'_A$ satisfies Condition~$(\star)$.

\item[(4)] $C,D\in\alpha$; \ $E\in\alpha_1\cup\alpha_2$.

Set $C':=\proj_AC$ and $D':=\proj_AD$. As in (3), we have $C',D'\in\alpha_1$. Then, again using the above
mentioned standard property of projections, we obtain

$$\begin{array}{rclr}
\delta(C,D)\delta(D,E)&=&\delta(C,C')\delta(C',D)\delta(D,E) & \hspace{1cm} \mbox{since }C,C',D\in\alpha,\\
                      &=&\delta(C,C')\delta(C',D)\delta(D,D')\delta(D',E) & \mbox{by }(3),\\
                      &=&\delta(C,C')\delta(C',D')\delta(D',E) & \mbox{since }C',D,D'\in\alpha,\\
                      &=&\delta(C,C')\delta(C',E) & \mbox{since }C',D',E\in\alpha_1\cup\alpha_2,\\
                      &=&\delta(C,E) & \mbox{by }(3).
\end{array}$$
\end{itemize}

Now the other cases all follow easily:

\begin{itemize}
\item[(5)] $C\in\alpha_1\cup\alpha_2$; \ $D,E\in\alpha$. The claim follows from (4) by taking inverses.

\item[(6)] $C,E\in\alpha$; \ $D\in\alpha_2$. By (4), we have $\delta(C,D)=\delta(C,E)\delta(E,D)$, hence
$$\delta(C,D)\delta(D,E)=\delta(C,E)\delta(E,D)\delta(D,E)=\delta(C,E).$$

\item[(7)] $C,D\in\alpha_1\cup\alpha_2$; \ $E\in\alpha$. Follows from (3) by taking inverses.

\item[(8)] $C,E\in\alpha_1\cup\alpha_2$; \ $D\in\alpha$. Using (3), we have
$$\delta(C,D)\delta(D,E)=\delta(C,D)\delta(D,C)\delta(C,E)=\delta(C,E).$$
\end{itemize}

So, by Proposition~\ref{fact2}, there is an apartment $\Sigma$ of $\Delta$ with $\alpha\cup\alpha_2\subseteq\mathrm{Ch}(\Sigma)$. Hence
$\Sigma$ contains $M$ (the boundary of $\widehat{\alpha}$) and also $\Sigma_A$ (since
$\cC'_A\subseteq\mathrm{Ch}(\Sigma)$). Therefore, $\Sigma\cap L_A$ contains $\Sigma_A$. However, since
$\Sigma\cap L_A$ is also an apartment of $L_A$ (because $A\in\Sigma$), we must have $\Sigma\cap
L_A=\Sigma_A$, and the assertion follows.

STEP II. We now reduce the general case to the case handled in Step I.

We first observe that $\widehat{\alpha}\cap L_A$ is a root in $L_A$ with $M\cap L_A$ as bounding wall. Recall
that this is also the common wall of the roots $\widehat{\alpha}_1,\widehat{\alpha}_2$ in $L_A$ corresponding
to $\alpha_1,\alpha_2$, respectively. Let $C$ be the chamber of $\alpha$ which contains the panel $P$
introduced in the beginning of this proof. Without loss of generality, we may assume that $C\neq C_1$
(otherwise we interchange the roles of $\alpha_1$ and $\alpha_2$). Then,
by Lemma~\ref{lemma3}, the union of the two roots $\widehat{\alpha}\cap L_A$ and $\widehat{\alpha}_1$ is an
apartment $\widetilde{\Sigma}_A$ of $L_A$.

We can now apply Step I to $\alpha$ and to this apartment $\widetilde{\Sigma}_A$. Thus we obtain an apartment
$\widetilde{\Sigma}$ with $\alpha\cup\alpha_1\subseteq\mathrm{Ch}(\widetilde{\Sigma})$. Let
$\widetilde{\alpha}$ be the root opposite $\widehat{\alpha}$ in $\widetilde{\Sigma}$. Then $M$ is also the boundary wall of $\widetilde{\alpha}$.
Note that $\widehat{\alpha}_1$ is a root in $\widetilde{\Sigma}\cap L_A$, and that $\widehat{\alpha}_1$ shares the boundary
$M\cap L_A$ with the two opposite roots $\widetilde{\alpha}\cap L_A$ and $\widehat{\alpha}\cap L_A$. Therefore, $\widehat{\alpha}_1=\widetilde{\alpha}\cap L_A$
or $\widehat{\alpha}_1=\widehat{\alpha}\cap L_A$. Since $C_1\neq C$, it follows that $\widehat{\alpha}_1=\widetilde{\alpha}\cap L_A$, and hence
$\alpha_1\subseteq \mathrm{Ch}(\widetilde{\alpha})$.
Now we can apply Step I to $\mathrm{Ch}(\widetilde{\alpha})$ and $\Sigma_A$, and the lemma is proved. \qed

\section{Complete systems of apartments}\label{complete}

In this section, we prove our First Main Result. We denote by $\Delta$ an arbitrary thick building, and $\cA$ is the complete apartment system of $\Delta$.

The proof of our First Main Result will be by an induction the first step of which is in fact
Proposition~\ref{abr}. Hence we first proof Proposition~\ref{abr}, which we restate here. Note that this
proposition is contained in \cite{Abr-Bro:08} in the form of the exercises 5.83 and 5.84.

\begin{prop}\label{abrbis}
If $\kappa$ is a convex chamber subcomplex of an apartment $\Sigma$ of $\Delta$, then there exists an apartment $\Sigma'$ of
$\Delta$ such that $\kappa=\Sigma\cap\Sigma'$.
\end{prop}

\proof We may assume $\kappa\neq\Sigma$ (otherwise set $\Sigma'=\Sigma$). Let $\cM$ be the set of walls of
$\Sigma$ containing some panel $P\in\kappa$ with the property that $\kappa$ contains exactly one chamber
through $P$. The set $\cM$ is not empty since $\Sigma$ is connected and $\kappa\neq\Sigma$. We choose an index
set $J$ such that $\cM=\{M_j\st j\in J\}$ and, for each $j\in J$, a panel $P_j$ contained in $M_j\cap\kappa$
so that $P_j$ is contained in precisely one chamber of $\kappa$. Furthermore, we denote by $C_j,D_j$, $j\in J$,
the two chambers of $\Sigma$ containing $P_j$, with $C_j\in\kappa$ and consequently $D_j\notin\kappa$. Let
$\alpha_j$, $j\in J$, be the root of $\Sigma$ containing $C_j$ but not $D_j$,
which implies that $M_j$ is the boundary of $\alpha_j$. Since $\kappa$ is a convex
chamber subcomplex of the apartment $\Sigma$, we easily see that $\kappa\subseteq\alpha_j$, for all $j\in J$.
Hence $$\kappa\subseteq\bigcap_{j\in J}\alpha_j.$$ Now suppose, by way of contradiction, that some simplex $X$
is contained in $\alpha_j$, for all $j\in J$, but not in $\kappa$. By considering a minimal gallery joining $X$
to $\kappa$, and taking into account that the intersection of all $\alpha_j$, $j\in J$, is convex itself, we
see that we may assume that $X$ is a chamber adjacent to some chamber $X'$ of $\kappa$. By construction, there
exists $i\in J$ such that $X\cap X' \subseteq M_i$. We now obtain the contradiction $X,X'\in \alpha_i$.

Hence we have shown that $\kappa$ is the intersection of all roots $\alpha_j$, for $j$ ranging over $J$. %Now,
%for all $i\in J$, we define $\kappa_i$ as the intersection of all roots $\alpha_j$, for $j$ ranging over
%$J\setminus\{i\}$. Since $P_i\not\subseteq M_j$ for $i\neq j$, we have $D_i\in \alpha_j$. But clearly
%$D_i\notin \alpha_i$. So we have shown that $\kappa\subseteq\kappa_i\neq \kappa$, for all $i\in J$.
Since $\Delta$ is thick, we can choose, for all $j\in J$, a chamber $D_j'\supseteq P_j$, with $C_j\neq D_j'\neq
D_j$, and we consider the sets

$$L=\kappa\cup\{D_j\st j\in J\}, \mbox{ and }L'=\kappa\cup\{D_j'\st j\in J\}.$$

We now claim that the mapping $\varphi:\cC(L)\rightarrow\cC(L')$, which is the identity on $\cC(\kappa)$ and
maps $D_j$ onto $D_j'$, for all $j\in J$, is a $W$-isometry. To that aim, we consider for each pair $i,j\in
J$, $i\neq j$, a minimal gallery $\gamma_{ij}$ joining $C_i$ with $C_j$ (and entirely contained in $\kappa$).
We extend this to a gallery $$\lambda_{ij}=(D_i,\underbrace{C_i,\ldots,C_j}_{\gamma_{ij}},D_j).$$

Since $M_i\neq M_j$, and since neither $M_i$ nor $M_j$ cross the gallery $\gamma_{ij}$ (because
$\kappa$ is contained in $\alpha_i$ and $\alpha_j$), we see that $\lambda_{ij}$ is a gallery in $\Sigma$ crossing every wall of $\Sigma$ at most once.
This implies that $\lambda_{ij}$ is a minimal gallery. If we put $s_i=\delta(C_i,D_i)$, for all $i\in J$, then
$\delta(D_i,D_j)=s_i\delta(C_i,C_j)s_j$, for all pairs $i,j\in J$, $i\neq j$.

Now the gallery $$\lambda'_{ij}=(D'_i,\underbrace{C_i,\ldots,C_j}_{\gamma_{ij}},D'_j)$$ has the same type as
$\lambda_{ij}$ and hence is also reduced. Hence $\delta(D'_i,D'_j)=s_i\delta(C_i,C_j)s_j=\delta(D_i,D_j)$. It
now follows rather easily that $\varphi$ is a $W$-isometry (the proof of $\delta(D_i,X)=\delta(D_i',X)$,
for all $i\in I$ and all $X\in\mathrm{Ch}(\kappa)$, is similar to that of $\delta(D_i,D_j)=\delta(D_i',D_j')$ above).
Since $L$ is contained in an apartment, namely in
$\Sigma$, the extension theorem, or equivalently, Proposition~\ref{fact2}, implies that $L'$ is contained in
some apartment $\Sigma'\in\cA$. Define $\kappa':=\Sigma\cap\Sigma'$. Then, by the foregoing, $\kappa\subseteq
L\cap L'\subseteq\Sigma\cap\Sigma'=\kappa'$. Hence $\kappa'$ is a convex subcomplex of $\Sigma$ containing
$\kappa$, in particular it contains $C_i$, for all $i\in J$. But it does not contain $D_i$, because $D_i\notin\Sigma'$.
Consequently $\kappa'\subseteq \alpha_i$, for all $i\in J$. Hence $$\kappa'\subseteq\bigcap_{i\in
J}\alpha_i=\kappa$$ and so $\kappa'=\kappa$.

The proposition is proved. \qed

From now on, for the rest of this section, our standing hypothesis is that each panel of
$\Delta$ is contained in at least $4$ chambers.

\begin{mlemma}\label{lemma2}
For every apartment $\Sigma\in\cA$ and every wall $M$ in $\Delta$, there exists an
apartment $\Sigma'\in\cA$ containing $M$ and satisfying $\Sigma'\cap \Sigma=M\cap\Sigma$.
\end{mlemma}

\proof The proof has two main steps: the case $M\cap\Sigma=\emptyset$ is treated separately, and then the
general case is reduced to the first step using Lemma~\ref{lemma1}.

\textbf{Case I: $M\cap\Sigma=\emptyset$.}

We choose an arbitrary panel $P$ in $M$ and four distinct chambers $C_1,C_2,C_3,C_4$ containing $P$. Then, by
Corollary~\ref{cor}, there exist four roots $\alpha_1,\alpha_2,\alpha_3,\alpha_4$ with boundary wall $M$
containing $C_1,C_2,C_3,C_4$, respectively, such that $\Sigma_{ij}:=\alpha_i\cup\alpha_j$ is an apartment, for
all $i,j\in\{1,2,3,4\}$, $i\neq j$.

The apartment $\Sigma$ meets every root $\alpha_i$, $i\in\{1,2,3,4\}$, in a convex subcomplex. If there are two
distinct $k,\ell\in\{1,2,3,4\}$ such that $\Sigma\cap\alpha_k=\Sigma\cap\alpha_\ell=\emptyset$,
then the assertion follows by setting $\Sigma'=\Sigma_{k\ell}$. Suppose now that for some $j\in\{1,2,3,4\}$ the
intersection $\Sigma\cap\alpha_j$ has dimension at least $1$, and that for some
$i\in\{1,2,3,4\}\setminus\{j\}$, the intersection $\Sigma\cap\alpha_i$ is nonempty.  Consider the convex
subcomplex $\Theta:=\Sigma_{ij}\cap\Sigma$ of $\Sigma_{ij}$ of dimension at least $1$. Then $\Theta$ is a
chamber complex in its own right (see Proposition~1$(ii)$ in~\cite{Abr:94}), and, in particular, $\Theta$ is a
connected simplicial complex (since $\dim\Theta>0$). By assumption there exists a vertex $x$ in
$\Sigma\cap\alpha_i$, and a vertex $y$ in $\Sigma\cap\alpha_j$. So, we find a path $x=x_0,x_1,\ldots,x_n=y$ in
$\Theta$ connecting $x$ and $y$ (all $\{x_\ell,x_{\ell+1}\}$ are edges in $\Theta$). Then at least one of the
vertices $x_k$ is in $M$. Indeed, if $y\in\alpha_i$, then $y=x_n\in\alpha_i\cap\alpha_j=M$. If $y$ is not in
$\alpha_i$, then there exists a $k$ with $x_k\in\alpha_i$ and $x_{k+1}\notin\alpha_i$. In particular, $x_{k+1}$
is not in $M$, and hence it is in the
interior of $\alpha_j$, implying $\{x_k,x_{k+1}\}\in\alpha_j$. Hence $x_k\in\alpha_i\cap\alpha_j=M$. However,
we assumed $M\cap\Sigma=\emptyset$. So $M\cap\Theta\neq\emptyset$ is not possible, and we obtain a
contradiction.

So we may assume that for every $i\in\{1,2,3,4\}$, $\Sigma\cap\alpha_i$ is either $0$-dimensional or empty.
We may also assume that $\Sigma\cap\alpha_1$ is $0$-dimensional. Let $x$ be a vertex in $\Sigma\cap\alpha_1$. Then, for all
$i\in\{2,3,4\}$, $\Sigma\cap\Sigma_{1i}$, being a $0$-dimensional convex subcomplex of $\Sigma$,  is contained in the support of $x$ in $\Sigma$
by Proposition~\ref{fact1}. Now, by Proposition~1$(iii)$ of~\cite{Abr:94}, this
support, which is also a $0$-dimensional convex subcomplex of $\Sigma$, can have at most two vertices. Hence there
is at most one $i\in\{2,3,4\}$ with $\Sigma\cap\alpha_i\neq\emptyset$ since $\Sigma\cap\alpha_i$ is contained
in the support of $x$ in $\Sigma$ and the intersections $\Sigma\cap\alpha_j$, for $j\in\{1,2,3,4\}$, are pairwise disjoint
(because $M\cap\Sigma=\emptyset$). Hence we can put $\Sigma'=\alpha_k\cup\alpha_\ell$, with
$k$ and $\ell$ distinct in $\{2,3,4\}$ such that $\Sigma\cap\alpha_k=\Sigma\cap\alpha_\ell=\emptyset$.

\textbf{Case II: $M\cap\Sigma\neq\emptyset$.}

Choose a maximal simplex $A$ in $M\cap\Sigma$ and set $L_A:=\mathrm{Lk}_\Delta(A)$, $M_A:=M\cap L_A$ and
$\Sigma_A:=\mathrm{Lk}_\Sigma(A)=\Sigma\cap L_A$. Then $M_A$ is a wall in $L_A$, $\Sigma_A$ is an apartment of
$L_A$ and $M_A\cap\Sigma_A=\emptyset$ since $A$ is maximal in $M\cap\Sigma$. Applying Case I to $M_A$ and
$\Sigma_A$, we find an apartment $\Sigma'_A$ of $L_A$ containing $M_A$ and satisfying
$\Sigma'_A\cap\Sigma_A=\emptyset$. In view of Lemma~\ref{lemma1} we can find $\Sigma'\in\cA$ with
$\Sigma'_A=\Sigma'\cap L_A$ and $M\subseteq \Sigma'$. This implies that $A$ is maximal in $\Sigma'\cap \Sigma$
(noting $A\subseteq B\in \Sigma'\cap\Sigma$ implies $B\setminus A\in \Sigma'_A\cap\Sigma_A=\emptyset$). And
$\Sigma'\cap\Sigma$ is a convex subcomplex of $\Sigma'$. Then $\Sigma'\cap\Sigma\subseteq M$ by
Proposition~\ref{fact1}, because $A\in M$. Hence $\Sigma'\cap\Sigma=\Sigma'\cap\Sigma\cap M=(\Sigma'\cap
M)\cap\Sigma=M\cap\Sigma$. \qed

\begin{theorem}\label{theo1}
Let $\Delta$ be a building with the property that each of its panels is contained in at least four chambers, 
and let $\kappa$ be a convex subcomplex of an apartment $\Sigma$ of $\Delta$. Then there exists an apartment $\Sigma'$ of $\Delta$
such that $\Sigma'\cap\Sigma=\kappa$.
\end{theorem}

\proof The proof goes by induction on the codimension codim$(\kappa)$. The case codim$(\kappa)=0$ is settled by
Proposition~\ref{abrbis}. So let us assume codim$(\kappa)>0$. Choose a maximal simplex $A$ of $\kappa$. Recall
that $\kappa$ is a chamber complex (see \cite{Abr:94}, Proposition 1$(ii)$). This implies in particular $\dim A
= \dim\kappa<\dim\Sigma$. So there exists a panel $P$ in $\Sigma$ with $A\subseteq P$. Denote by $M$ the wall
of $\Sigma$ containing $P$. Then $\kappa\subseteq M$ by Proposition~\ref{fact1}.

Let $C_1,C_2$ be the two chambers of $\Sigma$ which contain $P$. Let $x_1$ be the vertex in $C_1\setminus P$,
i.e.~$C_1=P\cup\{x_1\}$, and set $A_1:=A\cup\{x_1\}$. Note that $x_1\notin M$ and hence $A_1\notin M$. Now let
$\kappa_1$ be the convex hull of $\kappa\cup\{A_1\}$ in $\Sigma$. Then $\dim\kappa_1=\dim A_1$. Indeed, let
$M'$ be any wall of $\Sigma$ containing $A_1$. Then $M'$ contains $\kappa$ by Proposition~\ref{fact1} and hence
also $\kappa_1$ since $M'$ is convex and $\kappa\cup\{A_1\}\subseteq M'$. Since this is true for any wall $M'$
of $\Sigma$ containing $A_1$, we see that $\kappa_1$ is contained in the intersection supp$A_1$ of all these
walls (and which we have called the \emph{support} of $A_1$ in $\Sigma$ above).  We now claim that $A_1$
is maximal in supp$A_1$. Indeed, if $A_1'$ is any element of $\Sigma$ properly containing $A_1$, we choose a
chamber $C'$ of $\Sigma$ containing $A_1'$ and a panel $P'$ of $C'$ with $A_1\subseteq P'$ and $A_1'\not\subseteq P'$.
Then the wall through $P'$ contains $A_1$ but not $A_1'$, since it does not contain $C'=P'\cup\{A_1'\}$. Hence the claim.
Therefore $\dim\kappa_1=\dim(\mathrm{supp}A_1)=\dim A_1 = 1+ \dim A$.

In particular, we can apply the induction hypothesis to $\kappa_1$, which gives us an apartment $\Sigma_1$ of
$\Delta$ satisfying $\Sigma_1\cap\Sigma=\kappa_1$.

Now let $D_1$ be a chamber of $\Sigma_1$ containing $A_1$, and let $P_1$ be the panel
$P_1=D_1\setminus\{x_1\}$. So $P_1$ contains $A$ but not $A_1$. Denote by $M_1$ the wall of $\Sigma_1$
containing $P_1$. We see that $M_1$ also contains $A$ but not $A_1$. By Proposition~\ref{fact1}, this implies
first of all $\kappa\subseteq M_1$ (note that $\kappa$ is also a convex subcomplex of $\Sigma_1$ since $\kappa$
is convex in $\Delta$). Secondly, $A$ is maximal in $M_1\cap\kappa_1$. Indeed, if not, then $M_1$ would contain
some simplex strictly bigger than $A$ and hence of dimension $\dim A +1=\dim\kappa_1$. But then
Proposition~\ref{fact1} would imply $\kappa_1\subseteq M_1$, contradicting $A_1\notin M_1$. Hence
$\dim(M_1\cap\kappa_1)=\dim A=\dim\kappa$.

We now claim that $M_1\cap\kappa_1=\kappa$. The inclusion $\kappa\subseteq M_1\cap\kappa_1$ is clear (we
verified $\kappa\subseteq M_1$ above). Now let $\beta$ be any root of $\Sigma_1$ containing $\kappa$. If
$\beta$ contains $x_1$, it also contains $A_1=A\cup\{x_1\}$ and hence $\kappa_1$. If $\beta$ does not contain
$x_1$, then $A$, which is in $\beta$ and which is joinable to $x_1$, must be contained in the wall $N$ which
bounds $\beta$. However, we showed above that $A$ is maximal in $M_1\cap\kappa_1$. Hence $A\in N$ implies
$M_1\cap\kappa_1\subseteq N$ by Proposition~\ref{fact1}. So for \emph{any} root $\beta$ of $\Sigma_1$
which contains $\kappa$ we have $M_1\cap\kappa_1\subseteq\beta$. Hence $M_1\cap\kappa_1$ is contained in the
intersection of all these roots $\beta$, and this intersection is equal to $\kappa$ because $\kappa$ is a
convex subcomplex of $\Sigma_1$. So we have also proved $M_1\cap \kappa_1\subseteq\kappa$, and the claim
follows.

Recall that $\kappa_1=\Sigma_1\cap\Sigma$. Hence the previous claim immediately implies
$M_1\cap\Sigma=M_1\cap\Sigma_1\cap\Sigma=M_1\cap\kappa_1=\kappa$. Applying Lemma~\ref{lemma2}, we find an
apartment $\Sigma'$ of $\Delta$ which contains $M_1$ and satisfies $\Sigma'\cap\Sigma=M_1\cap\Sigma=\kappa$.
\qed

\begin{cor}
If $\kappa$ is a convex subcomplex of $\Sigma$ with $\dim\kappa<\dim\Sigma$, then there exists a wall $M$ in
$\Sigma$ containing $\kappa$, and for each such wall $M$ there exists some wall $M'$ in $\Delta$ such that
$\kappa=M\cap M'$.
\end{cor}

\proof Choose the apartment $\Sigma'$ as in Theorem~\ref{theo1} such that $\Sigma'\cap\Sigma=\kappa$. Let $A$
be a maximal simplex of $\kappa$ and choose walls $M,M'$ in $\Sigma,\Sigma'$, respectively, with $A\in M$ and
$A\in M'$ (this is possible since $\dim A=\dim\kappa<\dim\Sigma=\dim\Sigma'$). Then by Proposition~\ref{fact1},
$\kappa\subseteq M$ and $\kappa\subseteq M'$. So we have $$\kappa\subseteq M\cap M'\subseteq
\Sigma\cap\Sigma'=\kappa,$$ and the assertion follows.\qed

\section{Incomplete systems of apartments}\label{incomplete}

In this section we show our Second Main Result.

We start with a proposition that motivates our restriction to consider solely \emph{finite} convex subcomplexes in
the sequel.

\begin{prop}\label{infinite}
Let $\Delta$ be a thick building and let $\Sigma$ be an apartment (in the complete apartment system $\cA$ of
$\Delta$). Let $\kappa$ be an infinite subset of $\Sigma$. Then
the family $\cA^*=\{\Sigma^*\in\cA\st \kappa\not\subseteq\Sigma^*\}$ is a system of apartments for $\Delta$. In
particular, $\cA'=\cA^*\cup\{\Sigma\}$ is also a system of apartments for $\Delta$, and $\kappa$ is contained in
a unique member of $\cA'$. Hence if $\kappa$ is any infinite proper convex subcomplex of $\Sigma$, it cannot 
be the intersection of all apartments of $\cA'$ in which it is contained.
\end{prop}

\proof Let $C,D$ be two chambers of $\Delta$. Then the convex hull $H$ of $C,D$ is contained in some apartment $\Sigma_H\in\cA$.
Since $H$ is convex, Proposition~\ref{abrbis} implies the existence of an apartment $\Sigma'_H\in\cA$ such that
$\Sigma_H\cap\Sigma'_H=H$. Since $H$ is finite, we have $\kappa\not\subseteq H$ and so $\kappa$ cannot be
contained in both of $\Sigma_H$ and $\Sigma_H'$. Hence at least one of these is a member of $\cA^*$.

The assertion now follows easily. \qed

\begin{theorem}\label{theo2}
Let $\Delta$ be a thick building of type $(W,S)$. Then the following are equivalent.
\begin{itemize}
\item[\rm(I)] For any apartment system $\cA$ of $\Delta$, the convex closure of two given chambers $C$ and $D$
is precisely the intersection of two apartments in $\cA$ containing $C$ and $D$.

\item[\rm(II)] For any apartment system $\cA$ of $\Delta$, every finite convex chamber subcomplex $\kappa$
contained in some member of $\cA$ is precisely the intersection of two apartments in $\cA$ containing $\kappa$.

\item[\rm(III)] For every finite subset $F$ of chambers of the standard Coxeter complex $\Sigma(W,S)$, there exist
two chambers $C,D\in\Sigma(W,S)$ such that $F$ is contained in the convex closure of $C$ and $D$.

\item[\rm(IV)] For every triplet $\{X,Y,Z\}$ of chambers of the standard Coxeter complex $\Sigma(W,S)$, there exist
two chambers $C,D\in\Sigma(W,S)$ such that $\{X,Y,Z\}$ is contained in the convex closure of $C$ and $D$.
\end{itemize}

\end{theorem}

\proof We show (IV)$\Rightarrow$(III)$\Rightarrow$(II)$\Rightarrow$(I)$\Rightarrow$(IV), where the last implication 
requires most of the work.

(1) (IV)$\Rightarrow$(III). We use induction on $|F|\geq 3$. For $|F|=3$, this is precisely (IV). Now let
$|F|>3$ and choose $Z\in F$ arbitrary. Then, by the induction hypothesis, there exist two chambers $X,Y$ in
$\Sigma(W,S)$ such that $F\setminus\{Z\}$ is contained in the convex closure of $X$ and $Y$. Applying (IV) to
$X,Y,Z$ gives us the desired pair $C,D$ of chambers.

(2) (III)$\Rightarrow$(II). Denote by $\widetilde{\cA}$ the complete apartment system of $\Delta$. Suppose
$\kappa\subseteq\Sigma\in\cA$. By Proposition~\ref{abrbis}, there exists a member $\widetilde{\Sigma}$ of
$\widetilde{\cA}$ such that $\kappa=\Sigma\cap\widetilde{\Sigma}$. We may assume $\widetilde\Sigma\notin\cA$.
Define $$\widetilde\kappa:=\{E\in\mathrm{Ch}(\widetilde\Sigma)\st E \mbox{ has a panel in }\kappa\}.$$ Since
$\widetilde\kappa$ is clearly a finite set, (III) implies that $\widetilde\kappa$ is contained in the convex
hull $H$ of two chambers $X,Y$ of $\widetilde{\Sigma}$. Let now $\Sigma'$ be a member of $\cA$ containing $X,Y$. We
claim that $\kappa=\Sigma\cap\Sigma'$. Indeed, we clearly have $\kappa\subseteq\Sigma\cap\Sigma'$, so assume by
way of contradiction that $\Sigma\cap\Sigma'\not\subseteq\kappa$. Then there are adjacent chambers $C,D$ with
$C\in\kappa$, $D\notin\kappa$ and $C,D\in\Sigma\cap\Sigma'$. By definition of $\widetilde{\kappa}$, each panel $P$
contained in $\kappa$ is contained in precisely two chambers of $\widetilde{\kappa}$. Since $\widetilde{\kappa}
\subseteq H\subseteq \Sigma'$, these are the two chambers in $\Sigma'$ which contain $P$. So the two chambers of
$\Sigma'$ containing $P$ are in fact contained in $\widetilde{\kappa}$. Applied to $P=C\cap D$,
we obtain $C,D\in\widetilde{\kappa}$.
%Since $\widetilde\kappa\subseteq\Sigma'$, we have
%$$\widetilde\kappa=\{E\in\mathrm{Ch}(\Sigma')\st E \mbox{ has a panel in }\kappa\}.$$ Since $D$ has a panel in
%$\kappa$, it follows that $D\in\widetilde\kappa\subseteq\widetilde\Sigma$. hence
Hence $D\in\Sigma\cap\widetilde\Sigma=\kappa$, a contradiction. The claim is proved, and so is Assertion (II).

(3) (II)$\Rightarrow$(I). This is obvious since the convex closure of two chambers always contains a finite
number of chambers.

(4) (I)$\Rightarrow$(IV). Given three chambers $X,Y,Z$ in some apartment $\Sigma$ of $\Delta$, we show that
there exist two chambers $C,D$ in $\Sigma$ such that $X,Y,Z$ are contained in the convex hull of $C$ and $D$.
This is trivial if $Z$ is contained in the convex closure of $X$ and $Y$. So we may
assume that $Z$ is not contained in that convex closure. We first prove (IV) in the special case that $Z$ is
adjacent to some chamber $E$ of the convex closure of $X$ and $Y$. So assume, by way of contradiction, that
there are three chambers $C_1,C_2,C_3$ in some apartment $\Sigma$, with $C_3$ adjacent to some chamber $E$ of
the convex closure $\theta$ of $C_1$ and $C_2$, but $C_3$ does not belong to $\theta$,  such that the convex
closure of any two chambers $C,D$ of $\Sigma$ does not contain all of them. %Then, since the convex closure of
%two chambers is always contained in some apartment, and since apartments are isometric using an isometry that
%fixes pointwise their intersection, we see that for any pair $\{C,D\}$ of chambers of $\Delta$, the convex
%closure of $C$ and $D$ does not contain all of $C_1,C_2,C_3$.

Let $\widetilde{\cA}$ be the complete system of apartments for $\Delta$, and let $\cA$ be the subset of
$\widetilde{\cA}$ consisting of those apartments that either do not contain $\{C_1,C_2\}$ or contain
$\{C_1,C_2,C_3\}$. We show that $\cA$ is a system of apartments for $\Delta$. It suffices to show that any two
chambers $U,V$ are contained in a member of $\cA$. Let $\kappa$ be the convex closure of $U$ and $V$. Choose an
apartment $\Sigma'\in\widetilde{\cA}$ with $U,V\in\Sigma'$ and hence $\kappa\subseteq\Sigma'$. If
$\{C_1,C_2\}\not\subseteq\kappa$, then we choose, using Proposition~\ref{abrbis}, another apartment
$\Sigma''\in\widetilde{\cA}$ such that $\kappa=\Sigma'\cap\Sigma''$. Then not both of $\Sigma'$ and $\Sigma''$
contain $\{C_1,C_2\}$, and hence at least one of $\Sigma'$ and $\Sigma''$ is in $\cA$.

Suppose now that $\{C_1,C_2\}\subseteq\kappa$. Let $E'$ be the chamber of $\Sigma'$ which contains the panel
$P=E\cap C_3$ and is distinct from $E$. We first show that $E'\notin\kappa$. For this, let $\varphi:\Sigma'\rightarrow\Sigma$
be the isomorphism fixing $\Sigma'\cap\Sigma$ (and hence $\theta$) pointwise (see Proposition~\ref{newprop}).
Then the convex hull $\varphi(\kappa)$ of $\varphi(U)$ and $\varphi(V)$ contains $\varphi(\theta)=\theta$.
If $E'\in\kappa$, then $\varphi(\kappa)$ also contains $\varphi(E')$. Since $\varphi(E')$ is a chamber of $\Sigma$
with $\varphi(E')\neq\varphi(E)=E$ and $P=\varphi(P)\subseteq\varphi(E')$, we must have $\varphi(E')=C_3$. So the convex
closure of $\varphi(U)$ and $\varphi(V)$ contains $\varphi(\{C_1,C_2,E'\})=\{C_1,C_2,C_3\}$, contradicting our assumption.
Consequently we must have $E'\notin\kappa$.

Now let $\alpha$ be the root in $\Sigma'$ containing $E$ but
not $E'$. Since $\kappa$ is convex, $E\in\kappa$ and $E'\notin\kappa$, it follows that $\kappa\subseteq\alpha$.
Since we assume $\{C_1,C_2\}\subseteq\kappa$, also $\theta\subseteq\alpha$. By Proposition~\ref{extension}, there exists an
apartment $\Sigma''\in\widetilde{\cA}$ with $\alpha\cup\{C_3\}\subseteq\Sigma''$. Since $\{C_1,C_2,C_3\}\subseteq\Sigma''$,
we have $\Sigma''\in\cA$. This completes the proof that $\cA$ is a system of apartments for $\Delta$.

However, the convex closure of $C_1$ and $C_2$ can not be equal to the intersection of two
apartments of $\cA$ containing both of $C_1$ and $C_2$ since the former does not contain $C_3$ and the latter
always does. This contradicts Assumption~(I). So $C_1,C_2,C_3$ as described cannot exist. Consequently, if $Z$ is adjacent to a
chamber of the convex closure of $X$ and $Y$, then $\{X,Y,Z\}$ is contained in the convex closure of two
chambers of $\Sigma$.

We now prove the general case by induction on the gallery distance $n$ from $Z$ to the convex closure of $X,Y$. If
$n=1$, then this is the foregoing. If $n>1$, then there is a chamber $Z'$ adjacent to $Z$ and at distance
$n-1$ from the convex closure of $X,Y$. The induction hypothesis implies the existence of two chambers
$C',D'$ in $\Sigma$ such that $X,Y,Z'$ are contained in the convex closure of $C',D'$. Applying the above to
$\{C',D',Z\}$, the assertion is now clear.

Since every apartment of $\Delta$ is isomorphic to the standard Coxeter complex $\Sigma(W,S)$, the equivalence of (I) up to (IV) is
proved completely. \qed

For the classification of Coxeter complexes satisfying Condition (III) of the above theorem, we refer to
\cite{Cap:05}, Theorem 7.2. Necessary and sufficient is that each irreducible component is affine or
spherical. Since spherical buildings are automatically endowed with the complete system of apartments, we may
restrict to the irreducible affine case to state the following proposition.

\begin{prop}\label{cor2} Let $\Delta$ be an irreducible building of affine type with given apartment system $\cA$, and
suppose each panel is contained in at least four chambers. Then every finite convex subcomplex $\kappa$ contained in
some apartment $\Sigma\in\cA$ is the intersection of $\Sigma$ with some other member $\Sigma'\in\cA$.
\end{prop}

\proof  Let $\widetilde\cA$ be the complete system of apartments. By Theorem~\ref{theo1}, there exists
$\widetilde\Sigma \in\widetilde\cA$ such that $\kappa=\Sigma\cap\widetilde\Sigma$. Now denote by $\mathrm{Vert}(\kappa)$
the set of vertices of $\kappa$ and define
$$\widetilde\kappa:=\bigcup_{x\in\mathrm{Vert}(\kappa)}\mathrm{St}_{\widetilde\Sigma}(x).$$
Then $\widetilde\kappa$ is the union of a finite set of chambers of $\widetilde\Sigma$. Hence $\widetilde{\kappa}$ is
contained in the convex hull of two chambers, and there exists
$\Sigma'\in\cA $ with $\widetilde\kappa\subseteq\Sigma'$. For every $x\in\mathrm{Vert}(\kappa)$, the links
$\mathrm{Lk}_{\widetilde{\Sigma}}(x)$ and $\mathrm{Lk}_{{\Sigma'}}(x)$ are two apartments in $\mathrm{Lk}_\Delta(x)$ with
$\mathrm{Lk}_{\widetilde{\Sigma}}(x)\subseteq\mathrm{Lk}_{{\Sigma'}}(x)$, since $\mathrm{St}_{\widetilde{\Sigma}}(x)\subseteq\Sigma'$.
Hence $\mathrm{Lk}_{\widetilde{\Sigma}}(x)=\mathrm{Lk}_{{\Sigma'}}(x)$, and consequently also $\mathrm{St}_{\widetilde{\Sigma}}(x)=
\mathrm{St}_{{\Sigma'}}(x)$, implying
$$\widetilde\kappa=\bigcup_{x\in\mathrm{Vert}(\kappa)}\mathrm{St}_{\Sigma'}(x).$$ We claim that
$\Sigma\cap\Sigma'=\kappa$. Indeed, since $\kappa$ is certainly contained in both of $\Sigma$ and $\Sigma'$, we
may assume by way of contradiction that some vertex $y$ is contained in $\Sigma\cap\Sigma'$ but not in
$\kappa$.  By connectivity of $\Sigma\cap\Sigma'$ (recall that apartments in affine buildings are \emph{geodesically} convex),
this implies that there is some edge $\{x,x'\}$, with
$x\in\kappa$, and with $x'\in\Sigma\cap\Sigma'$ and not in $\kappa$.  Hence $x'\in\mathrm{St}_{\Sigma'}(x)$ and so
$x'\in\widetilde\kappa\subseteq\widetilde\Sigma$. Since $x'$ also
belongs to $\Sigma$, we obtain the contradiction $x'\in\Sigma\cap\widetilde\Sigma=\kappa$.

The proposition is proved. \qed

\begin{remark}\rm If $\kappa$ contains a chamber, the above proof yields (referring to Proposition~\ref{abrbis} instead of Theorem~\ref{theo1})
the same result without the assumption that each chamber of the thick building $\Delta$ is contained in at least four chambers.
\end{remark}

%\begin{remark}\rm The condition in Corollary~\ref{cor2} of $\Delta$ (which is not spherical)
%being of affine type is no longer a necessary condition if one restricts to finite convex subcomplexes of
%codimension at least one. Indeed, in this case all buildings of rank $3$ satisfy the conclusion of the
%corollary.
%\end{remark}

We now prove Proposition~\ref{propmain}.

\begin{prop}\label{adjacent}
For every thick building $\Delta$ and every apartment system $\cA$ of $\Delta$, the subcomplex formed by two adjacent
chambers is always the intersection of all members of $\cA$ containing both of those chambers.
\end{prop}

\proof
Let $C$ and $D$ be two adjacent chambers of $\Delta$. Let $E$ be a chamber in the intersection of all
apartments containing both of $C$ and $D$, and assume, by way of contradiction that $E\notin\{C,D\}$. By
convexity, we may assume that $E$ is adjacent to one of $C,D$, and without loss of generality, we may assume
$E$ is adjacent to $D$, and hence not to $C$. Let $E'$ be a third chamber containing the panel $D\cap E$,
$E'\notin\{D,E\}$. Then any apartment containing $C$ and $E'$ contains $D$ (since $(C,D,E')$ is a minimal gallery)
and does not contain $E$, a contradiction.

The proposition is proved completely. \qed

We end this paper with an application of the last result.

\begin{prop}
Let $\Delta$ be a thick building, and let $\Delta'$ be an arbitrary building. Let $\Delta$ and $\Delta'$ be endowed
with arbitrary systems of apartments $\cA$ and $\cA'$, respectively. Let $\varphi$ be a map from $\mathrm{Ch}(\Delta)$
to $\mathrm{Ch}(\Delta')$. Suppose that $\varphi$ bijectively maps the set
of chambers of any apartment $\Sigma\in\cA$ to the set of chambers of some apartment $\Sigma'\in\cA'$. Then $\varphi$ is injective
and preserves adjacency of chambers (i.e., any two chambers of $\Delta$ are adjacent if and only if their images under $\varphi$ are).
Furthermore, $\varphi(\mathrm{Ch}(\Delta))$ is the set of chambers of a thick subbuilding of $\Delta'$. Hence $\varphi$ induces a 
simplicial isomorphism of $\Delta$ onto a thick
subbuilding of $\Delta'$ in case $\Delta$ is $2$-spherical.
\end{prop}

\proof First note that $\varphi$ is injective. Indeed, this follows from the fact that $\varphi$ is injective
when restricted to the set of chambers of any apartment, and the fact that every pair of chambers is contained
in at least one apartment. Secondly, note that the image $\varphi(\mathrm{Ch}(\Delta))$ is the chamber set of a
subbuilding of $\Delta'$ with $\varphi(\cA)$ (obvious notation) as apartment system. Indeed, we only need to
check that every pair $\varphi(C),\varphi(D)$ of chambers of $\varphi(\mathrm{Ch}(\Delta))$ is contained in an
element of $\varphi(\cA)$ (with $C,D$ chambers of $\Delta$). But that follows immediately from our assumptions.
Hence, form now
on, we may assume that $\varphi(\mathrm{Ch}(\Delta))$ coincides with $\mathrm{Ch}(\Delta')$, and that $\cA'$ is
the set of all images under $\varphi$ of the members of $\cA$.

We now show that $\Delta'$ is thick.
Let $P$ be a panel of $\Delta'$ contained in the two chambers $\varphi(C)$ and $\varphi(D)$. Then, considering
any chamber of $\Delta$ adjacent to $C$ and different from $D$, we can use Proposition~\ref{adjacent} to see
that there exists some apartment $\Sigma\in\cA$ containing $C$ but not $D$. Let $E'$ be the unique chamber of
$\varphi(\mathrm{Ch}(\Sigma))$ adjacent to $\varphi(C)$ and containing $P$ (recall that $\varphi(\mathrm{Ch}(\Sigma))=
\mathrm{Ch}(\Sigma')$, for some $\Sigma'\in\cA'$). Since $D\notin\Sigma$, it follows
from global injectivity that $E'\neq\varphi(D)$. This shows thickness of
$\Delta'$.

Now we show that $\varphi$
preserves adjacency. Let $C,D$ be adjacent in $\Delta$. Then the intersection of all apartments in $\cA'$ containing $\varphi(C)$ and
$\varphi(D)$ contains every chamber of the convex hull $H$ of $\varphi(C)$ and $\varphi(D)$. Let $\varphi(E)$ be
a chamber in $H$. Then $E$ is contained in every apartment in $\cA$ containing $C,D$. By
Proposition~\ref{adjacent}, $E\in\{C,D\}$, and hence $\varphi(E)\in\{\varphi(C),\varphi(D)\}$. So we have proved that
$\varphi(C)$ and $\varphi(D)$ are the only chambers in $H$, showing that
$\varphi(C)$ and $\varphi(D)$ are adjacent.
Interchanging the roles of $\Delta$ and $\Delta'$, we also see that $\varphi^{-1}$ maps adjacent chambers to adjacent
chambers.

In view of Proposition 3.21 of \cite{Tit:74}, the proposition is completely proved. \qed

Addresses of the authors:

Peter Abramenko\\
Department of Mathematics\\
University of Virginia\\
P.O.Box 400137\\
Charlottesville, VA 22904\\
USA\\
\texttt{pa8e@virginia.edu}\\[15pt]

Hendrik Van Maldeghem\\
Department of Pure Mathematics and Computer Algebra\\
Ghent University\\
Krijgslaan 281, S22\\
9000 Gent\\
Belgium\\
\texttt{hvm@cage.UGent.be}


\begin{thebibliography}{99}

\bibitem{Abr:94} P. Abramenko, Walls in Coxeter complexes, \emph{Geom. Dedicata} \textbf{49} (1994), 71--84.

\bibitem{Abr-Bro:08} P.~Abramenko and K.~Brown, \emph{Buildings: Theory and Applications}, Graduate Texts in
Mathematics~\textbf{248}, Springer, New York, 2008.

\bibitem{Abr-Mal:04} P. Abramenko and H. Van Maldeghem, Maps between buildings that preserve a given Weyl distance,
\emph{Indag. Math.} \textbf{15} (2004), 305--319.

%\bibitem{Bro:89} K.~S.~Brown, \emph{Buildings}, Springer-Verlag New York, 1989.

\bibitem{Cap:05} P.-E.~Caprace, ``Abstract'' homomorphisms of split Kac-Moody groups. Preprint (December 2005), to appear as a Memoir of the AMS.

\bibitem{Ron:89} M.~A.~Ronan, \emph{Lectures on Buildings}, in: Perspect.~Math. \textbf{7},
Academic Press, San Diego, CA, 1989.

\bibitem{Tit:74} J. Tits, {\em Buildings of Spherical Type and Finite $BN$-pairs}, Lecture Notes in Math.\ {\bf 386},
Springer, Berlin--Heidelberg--New York, 1974.

\bibitem{Tit:81} J. Tits, A local approach to buildings, {\bf in} {\em The
Geometric Vein. The Coxeter Festschrift} (ed.\ D. Chandler {\em et al.}), Springer-Verlag (1981), 519--547.

\bibitem{Tit:86} J. Tits, Immeubles de type affine, \textbf{in} ``Buildings and the Geometry of Diagrams'', Como 1984
(ed. L. Rosati), \emph{Lect. Notes} \textbf{1181}, Springer-Verlag, Berlin, Heidelberg (1986), 157--190.

\bibitem{Tit:92} J. Tits, Twin buildings and groups of Kac-Moody type, \textbf{in} ``Groups, Combinatorics and
Geometry'', Durham 1990 (ed. M. Liebeck \& J. Saxl), \emph{London Math. Soc. Lecture Notes Ser.} \textbf{165},
Cambridge University Press, Cambridge (1992), 249--286.

\end{thebibliography}
\end{document}